\documentclass[a4paper,10pt,leqno]{article}
\usepackage[mathscr]{eucal}
\usepackage{amssymb}
\usepackage{latexsym}
\usepackage{theorem}
\theoremstyle{plain}

\newtheorem{Theorem}{Theorem}[section]
\newtheorem{Corollary}{Corollary}[section]
\newtheorem{Lemma}{Lemma}[section]
\newtheorem{Proposition}{Proposition}[section]

\theorembodyfont{\rmfamily}
\newtheorem{Remark}{Remark}[section]

\newtheorem{Definition}{Definition}[section]

\title{ 
Invariant minimal surfaces in 
the real special linear group of degree 2
}
\author{Jun-ichi Inoguchi
\footnote
{partially supported by Grand-in-Aid for
Encouragement of Young Scientists 12740051, 
14740053,
Japan Society for
Promotion of Science}}
\date
{{\it Dedicated to professor Koichi Ogiue
on his 60th birthday}
}
\begin{document}
\maketitle

\begin{abstract}
Invariant
minimal surfaces in 
the real special linear group
$\mathrm{SL}_{2}\mathbf{R}$ with canonical
Riemannian and Lorentzian metrics
are studied.

Constant mean curvature surfaces with 
vertically harmonic Gau{\ss} map are
classified.
\end{abstract}

\section*{Introduction}
In our previous works \cite{IKOS3}--\cite{IKOS4},
we have investigated fundamental properties of
the real special linear group
$\mathrm{SL}_{2}\mathbf{R}$ furnished
with canonical
left invariant Riemannian metric.
It is known that $\mathrm{SL}_{2}\mathbf{R}$
with canonical Riemannian metric admits
a structure of naturally reductive homogeneous
space and left invariant Sasaki structure.
The isometry group of the canonical 
left invariant metric is $4$-dimensional.

On the other hand, it is well known that
the Killing form of $\mathrm{SL}_{2}\mathbf{R}$ induces
a biinvariant Lorentz metric of constant curvature on 
$\mathrm{SL}_{2}\mathbf{R}$.

Thus $\mathrm{SL}_{2}\mathbf{R}$ with 
biinvariant metric is identified with
anti de Sitter $3$-space $H^3_1$.

As we will see in Section 1, 
the canonical left invariant Riemannian metric and
biinvariant Lorentzian metric (of constant
curvature $-1$) belong to same one-parameter family
of left invariant semi-Riemannian metrics.
Based on this fact,
in this paper, we shall give
a unified approach to geometry of 
$H^3_1$ and $\mathrm{SL}_{2}\mathbf{R}$
with canonical metric.

\vspace{0.2cm}

Since the canonical left invariant metric
is of non-constant curvature,
geometry of surfaces in $\mathrm{SL}_{2}\mathbf{R}$
is complicated.

In fact, 
we have shown in \cite{BDI1},
there are no extrinsic spheres
(totally umbilical surfaces with
constant mean curvature), especially
no totally geodesic surfaces in 
$\mathrm{SL}_{2}\mathbf{R}$.

\vspace{0.2cm}

In  \cite{Ko}, Kokubu introduced the notions of
{\it rotational surface} and {\it conoid} in  
$\mathrm{SL}_{2}\mathbf{R}$
with canonical left
invariant Riemannian metric.
Further he classified constant
mean curvature 
rotational surfaces and
minimal conoids.

Gorodski \cite{Grod} independently
investigated constant mean 
curvature rotational surfaces.

In \cite{BDI1},
Belkhelfa, Dillen and the 
author gave a characterisation of
rotational surfaces with constant
mean curvature. 
More precisely
a surface in $\mathrm{SL}_{2}\mathbf{R}$
is congruent to a 
rotational surface of constant
mean curvature if and only if its
second fundamental form is parallel.

\vspace{0.2cm}

In this paper we 
give some other characterisations of
rotational surfaces (of
constant mean curvature).

First we show that rotational 
surfaces in the sense of Kokubu coincide
with Hopf cylinders (over curves in the hyperbolic
$2$-space $H^2$)
in the sense of
Pinkall \cite{P} and 
Barros--Ferr{\'a}ndez--Lucas--Mero{\~n}o \cite{BFLM}.
Based on this fact,
we give a unified viewpoint for
\cite{BFLM} and \cite{Ko}.

Similarly we shall
show that conoids in the sense of Kokubu
coincide with Hopf cylinders over curves
in Lorentz $2$-sphere $S^2_1$.

\vspace{0.2cm}

When we identify the Lie algebra $\mathfrak{g}$
of $\mathrm{SL}_{2}\mathbf{R}$ with (semi) Euclidean 
$3$-space,
both $H^2$ and $S^2_1$ are given by
adjoint orbits in $\mathfrak{g}$.
The adjoint orbits of $\mathrm{SL}_2\mathbf{R}$
in $\mathfrak{g}$ are $H^2$, $S^2_1$ and lightcone
$\Lambda$. 
Based on this fact,
we shall introduce a new class of surfaces
in $\mathrm{SL}_{2}\mathbf{R}$.
More precisely, 
in section 4,
we shall investigate 
surfaces in $\mathrm{SL}_{2}\mathbf{R}$ derived from
curves in $\Lambda$. 

\vspace{0.2cm}

For every
surface in $\mathrm{SL}_2{\mathbf R}$,
we associate
a map into the Grassmannian bundle
$Gr_2(T\> \mathrm{SL}_{2}\mathbf{R})$
of $2$-planes---called
the {\it Gau{\ss}} {\it map} of the surface.
We shall give a characterisation
of constant mean curvature
rotational surfaces in terms of 
harmonicity for Gau{\ss} maps. 

More precisely, in the
final section, 
we shall prove that
a constant
mean curvature surface 
in $\mathrm{SL}_{2}\mathbf{R}$ is 
congruent to a rotational surface with
constant mean curvature if and only if
its Gau{\ss} map is vertically harmonic.

\vspace{0.15cm}

The author would like to thank 
professor
Luis Jose Al{\' \i}as (Universidad de Murcia,
Spain)
for his careful reading of the manuscript
and invaluable suggestions.

\section{The special linear group} 
{\bf 1.1} \hspace{0.15cm} Let 
$G=\mathrm{SL}_{2}{\mathbf R}$ be
the real special linear group of degree $2$:
$$
\mathrm{SL}_{2}\mathbf{R}=\left \{
\left (
\begin{array}{cc}
a & b \\ c & d
\end{array}
\right )
\ \biggr \vert \
a,b,c,d \in \mathbf{R},
\ ad-bc=1
\right \}.
$$

By using the Iwasawa decomposition 
$G=NAK$ of $G$; 
$$
N=\left \{
\left (
\begin{array}{cc}
1 & x \\ 0 & 1 
\end{array}
\right )
\
\biggr \vert \ 
x \in \mathbf{R}
\right \},
\eqno{\rm (Nilpotent \ part)}
$$
$$
A=
\left \{
\left (
\begin{array}{cc}
\sqrt{y} & 0 \\ 0 & 1/\sqrt{y}
\end{array}
\right )
\ \biggr
\vert \
y>0
\right \},
\eqno{\rm (Abelian \ part)}
$$
$$ 
K=
\left \{
\left (
\begin{array}{cc}
\cos \theta & \sin \theta \\
-\sin \theta & \cos \theta 
\end{array}
\right )
\ \biggr \vert 
\ 0\leq \theta \leq 2\pi
\right \},
\eqno{\rm (Maximal \ torus)}
$$
we can introduce the following
global coordinate system $(x,y,\theta)$
of $G$:
\begin{equation} \label{coord}
(x,y,\theta)\longmapsto
\left (
\begin{array}{cc}
1 & x \\ 0 & 1 
\end{array}
\right )
\left (
\begin{array}{cc}
\sqrt{y} & 0 \\ 0 & 1/\sqrt{y}
\end{array}
\right )
\left (
\begin{array}{cc}
\cos \theta & \sin \theta \\
-\sin \theta & \cos \theta 
\end{array}
\right ).
\end{equation}
We equip on $G$ the following 
one-parameter family
$\{g[\nu] \}$ of 
semi-Riemannian metrics:
$$
g[\nu]=\frac{dx^2+dy^2}{4y^2}+\nu
\left (
d \theta+\frac{dx}{2y}
\right )^2,\ \ \nu \in
\mathbf{R}^{*}.
$$
Every metric $g[\nu]$ is left invariant.
Clearly $g[\nu]$ is
Riemannian for $\nu>0$ 
and Lorentzian for $\nu<0$.

Throughout this paper we restrict our
attention to $\nu=\pm 1$
for simplicity.

One can see that $g[1]$ is only left invariant
but $g[-1]$ is a biinvariant Lorentz metric on $G$.

We take the following orthonormal coframe field
of $(G,g[\nu])$:
\begin{equation}\label{cobase}
\omega^1=\frac{dx}{2y},\ \
\omega^{2}=\frac{dy}{2y}, \ \
\omega^{3}=d\theta+\frac{dx}{2y}.
\end{equation}
The dual frame field 
of $\{\omega^1,\omega^2,\omega^3\}$ is given by
$$
\epsilon_1=2y
\frac{\partial}{\partial x}-
\frac{\partial}{\partial \theta},\
\epsilon_2=2y
\frac{\partial}{\partial y},\
\epsilon_3=
\frac{\partial}{\partial \theta}.
$$
Note that this orthonormal frame 
field is {\it not} left invariant.

The Levi-Civita connection $\nabla$ 
of $g[\nu]$ is given by
the following formulae:

$$
\nabla_{\epsilon_1}\epsilon_1=2\epsilon_2,\ \
\nabla_{\epsilon_1}\epsilon_2=-2\epsilon_1-\epsilon_3
,\ \
\nabla_{\epsilon_1}\epsilon_3=\nu\epsilon_2,
$$
\begin{equation}\label{con}
\nabla_{\epsilon_2}\epsilon_1=\epsilon_3,\ \
\nabla_{\epsilon_2}\epsilon_2=0
,\ \
\nabla_{\epsilon_2}\epsilon_3=-\nu\epsilon_1,
\end{equation}
$$
\nabla_{\epsilon_3}\epsilon_1=\nu\epsilon_2,\ \
\nabla_{\epsilon_3}\epsilon_2=-\nu \epsilon_1
,\ \
\nabla_{\epsilon_3}\epsilon_3=0.
$$
The commutation relations
of the basis are given by
\begin{equation}\label{Com}
[\epsilon_{1},\epsilon_{2}]=
-2\epsilon_{1}-2\epsilon_{3},\ \
[\epsilon_{1},\epsilon_{3}]=0,\ \
[\epsilon_{2},\epsilon_{3}]=0.
\end{equation}
The Riemannian curvature tensor
$R$ of the metric $g$
defined by
$$
R(X,Y)Z:=\nabla_{X}\nabla_{Y}Z-
\nabla_{Y}\nabla_{X}Z-
\nabla_{[X,Y]}Z,\ \ 
X,Y, Z \in \mathfrak{X}(G)
$$
is 
described by the
following formulae:

\begin{equation}\label{R}
\begin{array}{cc}
  {R}(\epsilon_{1},\epsilon_{2})
  \epsilon_{1}=(3\nu+4)\epsilon_{2}, &
 {R}(\epsilon_{1},\epsilon_{2})\epsilon_{2}
 =-(3\nu+4)\epsilon_{1}, \\
 {R}(\epsilon_{1},\epsilon_{3})\epsilon_{1}
 =-\nu\epsilon_{3}, &
 {R}(\epsilon_{1},\epsilon_{3})
 \epsilon_{3}=\nu \> \epsilon_{1}, \\
 {R}(\epsilon_{2},\epsilon_{3})
 \epsilon_{2}=-\nu \>\epsilon_{3}, &
{R}(\epsilon_{2},\epsilon_{3})
\epsilon_{3}=\nu \>\epsilon_{2}.
\end{array}
\end{equation}

\vspace{0.2cm}

\noindent
{\bf 1.2} \hspace{0.15cm} The
one-form $\eta=-d \theta-dx/(2y)$ 
is a {\it contact form} on
$G$, {\it i.e.},
$d \eta \wedge \eta \not=0$.

Let us define an endmorphism field $F$ by
$$
F\> \epsilon_1=\epsilon_2,\ 
F\> \epsilon_2=-\epsilon_1,\ 
F\> \epsilon_3=0.
$$
And put $\xi=-\epsilon_3$. Then
$(\eta,\xi,F,g[\nu])$ satisfies the 
following relations:
$$
F^2=-I+\eta \otimes \xi,
\ \
d \eta(X,Y)=2g(X,FY),
$$
$$
g(FX,FY)=g(X,Y)-\nu\eta(X)\eta(Y),
$$
$$
\nabla_{X}\xi=-\nu FX,
$$
$$
(\nabla_{X}F)Y=g(X,Y)\xi-\nu\eta(Y)X
$$
for all $X,Y \in \mathfrak{X}(G)$.

These formulae say that the
structure $(\xi,F,g[\nu])$ is
the associated almost contact
structure of the contact manifold
$(G,\eta)$ \cite{IKOS3}.
The resulting almost contact manifold
$(G;\eta,\xi,F,g[\nu])$ is a homogeneous
Sasaki manifold \cite{Tak}. 
The structure $(\eta,\xi,F,g[\nu])$ is called the
{\it canonical Sasaki structure} 
of $G$. With respect to the canonical
Sasaki structure, 
$(G,g[\nu])$ is a 
Sasaki manifold of constant holomorphic 
sectional curvature $-(3\nu+4)$.
The vector field $\xi$ is called the
{\it Reeb vector field} of $G$ associated to
$\eta$.
In Lorentzian case, 
since $\xi$ is a 
globally defined unit timelike
vector field on $G$,
$\xi$ time-orients $G$. 
\begin{Remark}
The Riemannian curvature tensor
$R$ of $(G,g[\nu])$
is given 
explicitly by 
\begin{eqnarray*}
R(X,Y)Z & = & -
g(Y,Z)X+g(Z,X)Y \\
& & -(1+\nu)\>\{
\eta(Z)\eta(X)Y
-\eta(Y)\eta(Z)X \\
& &+g(Z,X)\eta(Y)\xi
-g(Y,Z)\eta(X)\xi \\
& &-g(Y,FZ)FX-g(Z,FX)FY+
2g(X,FY)FZ
\>\}
\end{eqnarray*}
in terms of the canonical Sasaki structure.
In particular, this
explicit formula says
$g[-1]$ is a
Lorentz metric of 
constant curvature $-1$. 
As we will see later 
$(G,g[-1])$ is identified with the
anti de Sitter space $H^3_1$.

For more informations on 
the canonical Sasaki structure  
of $G$, we refer to \cite{IKOS3}.
\end{Remark}

\vspace{0.2cm}

\noindent
{\bf 1.3} \hspace{0.15cm} The 
special linear group $G$ 
acts transitively and
isometrically on the upper half plane:
$$
H^2(1/2)=
\left(
\{(x,y)\in {\mathbf R}^2 \ | \ y>0\},
\frac{dx^2+dy^2}{4y^2}
\right)
$$
of constant curvature $-4$.
The isotropy subgroup of $G$ at $(0,1)$ is 
the rotation group $K=\mathrm{SO}(2)$.
The natural projection $\pi:(G,g[\nu])
\to G/K=H^2(1/2)$
is a 
semi-Riemannian submersion
with totally geodesic fibres.
Moreover $\pi$ is given explicitly by
$$
\pi(x,y,\theta)=(x,y) \in H^2(1/2)
$$
in terms of the global coordinate
system (\ref{coord}).

The horizontal distribution 
of this semi-Riemannian 
submersion coincides
with the contact distribution 
determined by $\eta$.
The submersion 
$\pi:(G,g[-1]) \to H^2(1/2)$ 
is traditionally called the
{\it Hopf fibering} of $H^2(1/2)$. 
The Sasaki manifold
$(G,\eta;\xi,F,g[-1])$ is an example
of regular contact spacetime
which is {\it not} globally
hyperbolic.

\vspace{0.2cm}

\noindent
{\bf 1.4} \hspace{0.15cm}  Let
us denote by 
$\mathfrak{g}$ the Lie algebra of
$G$, {\it i.e.,} 
the tangent space of 
$G$ at the identity matrix
$\mathbf{1}$:
$$
\mathfrak{g}=
\left \{
\left (
\begin{array}{cc}
a & b \\ c & d
\end{array}
\right )
\ \biggr \vert \
a,b,c,d \in \mathbf{R},
\ a+d=0
\right \}.
$$

We take the following 
(split-quaternion) basis of
$\mathfrak{g}$:
$$
\mathbf{i}=
\left(
\begin{array}{cc}
0 & -1 \\
1 & 0
\end{array}
\right ),\ \
\mathbf{j}^\prime=
\left(
\begin{array}{cc}
0 & 1 \\
1 & 0
\end{array}
\right ),\ \
\mathbf{k}^\prime=
\left(
\begin{array}{cc}
-1 & 0 \\
0 & 1
\end{array}
\right ).
$$

Hereafter we identify 
$\mathfrak{g}$ with Cartesian $3$-space
$\mathbf{R}^3$ via the linear isomorphism:
$$
X=x_1\,{\mathbf i}+x_2\,{\mathbf j}^\prime
+x_3\,{\mathbf k}^\prime 
\longmapsto (x_1,x_2,x_3).
$$
Equivalently,
$$
X=\left(
\begin{array}{cc}
-x_3 & -x_1+x_2 \\
x_1+x_2 & x_3
\end{array}
\right )
\longmapsto
(x_1,x_2,x_3).
$$

We denote the scalar
product on ${\mathfrak g}$
induced by $g[1]$ and
$g[-1]$ by
$\langle \cdot,\cdot \rangle^{(+)}$
and 
$\langle \cdot,\cdot \rangle^{(-)}$
respectively.

The scalar products
$\langle \cdot,\cdot \rangle^{(\pm)}$ are
given explicitly by
the following formulae:
$$
\langle X,Y 
\rangle^{(+)}=
\frac{1}{2}
\mathrm{tr}
({}^{t}XY), 
\ X,Y \in \mathfrak{g},
$$
$$
\langle X,Y 
\rangle^{(-)}=
\frac{1}{2}\mathrm{tr}
(XY),
\ X,Y \in \mathfrak{g}.
$$
For $X\in \mathfrak{g}$,
$$
\langle X,X \rangle^{(\pm)}
=\pm x^2_1+x^2_2+x^2_3.
$$
Thus we identify 
$(\mathfrak{g},\langle
\cdot,\cdot
\rangle^{(+)})$
with Euclidean $3$-space:
$$
\mathbf{E}^3=
(\mathbf{R}^3(x_1,x_2,x_3),
dx_1^2+dx_2^2+dx_3^2).
$$
And 
$(\mathfrak{g},\langle
\cdot,\cdot \rangle^{(-)})$
is identified with 
Minkowski $3$-space: 
$$
\mathbf{E}^3_1=
(\mathbf{R}^3(x_1,x_2,x_3),
-dx_1^2+dx_2^2+dx_3^2)
$$
respectively.

Moreover the semi-Euclidean $4$-space
$$
\mathbf{E}^4_2=( \mathbf{R}^{4}(x_0,x_1,x_2,x_3) \ , \
-dx_0^2-dx_1^2+dx_2^2+dx_3^2\ )
$$ 
is identified with
the space $\mathrm{M}_{2}\mathbf{R}$
of all real $2$ by $2$ matrices:
$$
\mathrm{M}_{2}\mathbf{R}=\left
\{x_0 \mathbf{1}+x_1\mathbf{i}
+x_2 \mathbf{j}^{\prime}
+x_3 \mathbf{k}^{\prime}
\right\}.
$$
The semi-Euclidean metric of 
$\mathbf{E}^4_2$ corresponds to
the scalar product
$$
\langle X,Y 
\rangle=
\frac{1}{2}
\left\{
\mathrm{tr}
(XY)-\mathrm{tr}(X)
\mathrm{tr}(Y)
\right \},
\ X,Y \in \mathrm{M}_{2}\mathbf{R}.
$$
Since $\langle X,X
\rangle=-\det X$
for all $X 
\in \mathrm{M}_{2}\mathbf{R}$,
the special linear group
$G$ with biinvariant Lorentz 
metric $g[-1]$ is 
identified with
anti de Sitter $3$-space:
$$
H^3_1=\{ (x_0,x_1,x_2,x_3)
\in \mathbf{E}^4_2 \ |
\ -x_0^2-x_1^2+x_2^2+x_3^2=-1
\ \}.
$$
\vspace{0.2cm}

\noindent
{\bf 1.5} \hspace{0.15cm} The 
Lie group $G$ acts on $\mathfrak{g}$ by the 
Ad-action:
$$
\mathrm{Ad}:G \times \mathfrak{g}
\to \mathfrak{g};\
\mathrm{Ad}(a)X=a \> X \> a^{-1},\ a \in G,\ X \in
\mathfrak{g}.
$$

Since the determinant function $\det$ is Ad-invariant,
the Ad-orbits in $\mathfrak{g}$ are parametrised in 
the following way:

$$
\mathcal{O}_c=\left \{
X \in \mathfrak{g} \ \vert \
\det X=c
\ \right \},\ c \in \mathbf{R}.
$$
For $c\geq 0$, put
$$
\mathcal{O}^{\pm}_c=
\{ (x_1,x_2,x_3)  \in  \mathcal{O}_c
\vert \  \pm x_1>0
\}.
$$
Then 
\begin{eqnarray*}
\mathcal{O}_c &= &\mathcal{O}_c^{+}
\cup \mathcal{O}_c^{-},\ c>0, \\
\mathcal{O}_0&=&\mathcal{O}_0^{+}
\cup \{0\}
\cup \mathcal{O}_0^{-},\ c=0.
\end{eqnarray*}

\begin{Proposition}
The $\mathrm{Ad}$-orbits of $G$ are 
$$
\mathcal{O}^{\pm}_c,\ (c>0),\ \
\mathcal{O}_0^{\pm},\ (c=0), \
\{0\}, \ \
\mathrm{or}\ \ \mathcal{O}_c, \ (c<0).
$$

With respect to the Lorentz scalar product
$\langle \cdot,\cdot \rangle^{(-)}$,
the non-trivial $\mathrm{Ad}$-orbit $\mathcal{O}_c$ are classified as
follows{\rm :}

\vspace{0.2cm}

{\rm (1)} $c<0$:
The $\mathrm{Ad}$-orbit 
$\mathcal{O}_c$ is 
the pseudo-$2$-sphere
$S^2_1(\sqrt{-c})$
of radius $\sqrt{-c}$.

In this case $\mathcal{O}_c=G/A\mathbf{Z}^2$.

\vspace{0.2cm}

{\rm (2)} $c>0$:
The $\mathrm{Ad}$-orbit 
$\mathcal{O}^{\pm}_c$ is the upper or lower imbedding
of hyperbolic 
\par
$2$-space $H^2(\sqrt{c})$ with radius $\sqrt{c}$ 
in $\mathbf{E}^3_1$. 
In this case $\mathcal{O}^{\pm}_c=G/K$.

\vspace{0.2cm}

{\rm (3)} $c=0$:
The $\mathrm{Ad}$-orbit
$\mathcal{O}^{\pm}_0$ is the 
future or past lightcone{\rm :} 
$$
\Lambda_{\pm}=\{ (x_1,x_2,x_3)\not=0 \ \vert \
-x_1^2+x_2^2+x_3^2=0, \ \pm x_1>0
\}.
$$
The future lightcone $\Lambda_{+}$ 
is represented as
$\Lambda_{+}=G/N\mathbf{Z}_{2}$.
\end{Proposition}

\vspace{0.2cm}

\noindent
{\bf 1.6} \hspace{0.15cm} The Riemannian metric 
$g[1]$ is not only
$G$-left invariant but also right $K$-invariant.
Thus the product group
$G\times K$ acts isometrically on 
$(G,g[1])$.
Note that $(G,g[1])$
is represented by
$(G\times K/K,g[1])$
as a {\it naturally reductive}
(Riemannian) homogeneous space (See \cite{TV}).

On the other hand, 
since $g[-1]$ is biinvariant,
$G \times G$ acts 
isometrically on $(G,g[1])$.
Moreover $(G,g[-1])$ is represented by
$(G\times G/G,g[-1])$ as a
{\it Lorentzian symmetric} space.

Hence every subgroup of $G \times K$ 
acts isometrically on both
$(G,g[\nu])$.
Kokubu introduced the notion of
{\it helicoidal motion}
for $(G,g[1])$. This notion can be naturally
extended for $(G,g[\nu])$.

\begin{Definition} 
Let $\{\sigma^\mu_t\}_
{t \in {\mathbf R}}$ be a 
one parameter subgroup of $G\times K$
defined by
\begin{equation}
\sigma^\mu_t(X)=
\left(
\begin{array}{cc}
1 & \mu t \\ 0 & 1
\end{array}
\right)
X
\left(
\begin{array}{cc} 
\cos t & \sin t \\
-\sin t & \cos t
\end{array}
\right),\ \
\mu  \in \mathbf{R}.
\end{equation}
An element of $\{\sigma^\mu_t\}_{t \in {\mathbf R}}$ 
is called a {\it helicoidal motion with pitch} $\mu$.
\end{Definition}

Kokubu called surfaces in $(G,g[1])$ which are invariant under some
helicoidal motion group $\{\sigma^\mu_t\}$
{\it helicoidal surfaces}.

\section{Hopf cylinders}

\noindent
{\bf 3.1} \hspace{0.15cm}  
We recall two classes of surfaces in $(G,g[1])$
studied by Kokubu.

\begin{Definition}{\rm (\cite{Ko})}
An immersed surface in $G$ is said to be
a {\it rotational surface} if it is invariant under 
the right $K$-action.
\end{Definition}

A rotational surface can be parametrised as
\begin{equation}
\varphi(u,v)=
\left (
\begin{array}{cc}
1 & x(v) \\ 0 & 1 
\end{array}
\right )
\left (
\begin{array}{cc}
\sqrt{y(v)} & 0 \\ 0 & 1/\sqrt{y(v)}
\end{array}
\right )
\left (
\begin{array}{cc}
\cos u & \sin u \\
-\sin u & \cos u
\end{array}
\right).
\end{equation}
Obviously this definition is also
valid for $H^3_1$.

Next we recall the notion of
Hopf cylinder introduced by
Pinkall \cite{P}.

Let $\pi:S^3 \to S^2(1/2)$ be the Hopf fibering
of $S^2(1/2)$. 
Take a curve ${\bar \gamma}$ 
in the base space $S^2(1/2)$.
Then the inverse image 
$M:=\pi^{-1}\{{\bar \gamma}\}$
is a flat surface in $S^3$ which is called 
a {\it Hopf cylinder} over
${\bar \gamma}$ \cite{P}.
This construction is valid for other Hopf fiberings:
$
H^3_1 \to H^2(1/2)
$
and
$
H^3_1 \to S^2_1(1/2).
$

In particular Hopf  
cylinders in $H^3_1$ 
over curves in $H^2(1/2)$ are timelike.
Barros, Ferr{\'a}ndez, Lucas and
Mero{\~ n}o \cite{BFLM}, \cite{BFLM3},
\cite{F}
developed detailed studies on Hopf cylinders
in $H^3_1$.
It is easy to see that the notion of
Hopf cylinder can be extended naturally to
the fibering:
$$
\pi:(G,g[\nu]) \to H^2(1/2).
$$

 By using $\mathrm{SL}_{2}{\mathbf R}$-model 
of $H^3_1$ and the coordinate
system (\ref{coord}),
we can see that Hopf cylinders over
curves in $H^2$ are 
nothing but surfaces in $G$ 
invariant under the right action of $K$.

\begin{Proposition}
Let $M$ be a surface in $(G,g[\nu])$.
Then $M$ is a Hopf cylinder over a curve 
in $H^2(1/2)$ if and only if
it is a rotational surface.
\end{Proposition}

Thus we can unify two theories of ``Hopf cylinders in $H^3_1$" 
and of ``rotational surfaces in $(G,g[1])$".

\begin{Proposition}
Let $\varphi:I\times S^1
\rightarrow (G,g[\nu])$ be a Hopf cylinder
over a curve $(x(v),y(v))$ in $H^2(1/2)$
parametrised by arclength parameter $v$.
Then the induced metric of $\varphi$ is
\begin{equation}
\mathrm{I}[\nu]=
\nu \left( du+\frac{x^{\prime}(v)}{2y}dv
\right)^2+dv^2.
\end{equation}

Hence the Hopf cylinder $(I\times S^1,\varphi)$
is flat.
\end{Proposition}

Hopf cylinders of constant mean
curvature are classified as follows
(Compare \cite{BFLM}--\cite{BFLM3} and Proposition 4.3 in
\cite{Ko}):

\begin{Proposition}
{\rm (Classification 
of CMC Hopf cylinders)}

Let $c$ be a unit speed curve in $H^2(1/2)$ 
with curvature $\kappa$ and
$M_c$ the Hopf cylinder over $c$ in $(G,g[\nu])$.
Then $M_c$ is
of constant mean curvature if and only if
$c$ is a Riemannian circle in $H^2(1/2)$. 
The mean curvature of $M_c$ is
$H=\kappa/2$.

The Hopf cylinder $M_c$ is
classified in the 
following way{\rm :}
\newline {\rm(1)}
$M_c$ is a minimal complex circle if $\kappa=0$,
\newline {\rm(2)}
$M_c$ is a non-minimal complex circle 
or a Hopf cylinder over a line segment
\indent
$y=\pm (\sqrt{1-4\kappa^2}/(2\kappa))x$
if $0<\kappa^2<4$,
\newline {\rm(3)}
$M_c$ is a Hopf cylinder over a horocycle
or $y=$constant if $\kappa^2=4$,
\newline {\rm(4)}
$M_c$ is an imbedded torus
if $\kappa^2>4$.
\end{Proposition}
Note that,
in $H^3_1$ case, $M_c$ is a
$B$-scroll 
of the horizontal lift ${\hat c}$
of $c$. (\cite{DN}, \cite{BFLM}.
Compare with Theorem \ref{umbilicScroll}).

\begin{Remark}
The notion of complex circle is introduced by
Magid. (See \cite{Ma}, Example 1.12.) 
The non-minimal complex circle is an 
isometric immersion
$\varphi:{\mathbf E}_1^2(u,v)
\to H^3_1$ of Minkowski plane into
$H^3_1$ defined by
\[
\varphi(u,v)=
\left (
\begin{array}{c}
b\cosh v \cos u-a \sinh v \sin u \\
a\sinh v \cos u+b \cosh v \sin u \\
a\cosh v \cos u+b \sinh v \sin u \\
a\cosh v \sin u-b \sinh v \cos u
\end{array}
\right),
\]
where $a^2-b^2=-1,\ ab\not=0$.
The non-minimal complex circle $\varphi$ is a 
non-minimal flat timelike surface in $H^3_1$.
({\it cf}. Al{\' \i}as, Ferr{\'a}ndez and Lucas \cite{AFL},
Example 3.3.)

If we interchange $+$ and $-$ in the
third and fourth components of $\varphi$,
then we obtain a timelike minimal surface in $H^3_1$.
This timelike minimal surface has the following
expression $
\exp(u{\mathbf i})\
\exp(v{\mathbf k}^\prime)\
\exp(t{\mathbf j}^{\prime})$.
Here we put 
$b=\cosh t$ and $a=\sinh t$. 
\end{Remark}

\begin{Remark}\label{parallel}
It is straightforward to check that
every rotational surface of
constant mean curvature
in $(G,g[1])$
has parallel second fundamental form
(especially constant principal curvatures).
Conversely, one can see that surfaces with 
parallel second fundamental form in 
$(G,g[1])$ are congruent to 
rotational surfaces of 
constant mean curvature.
See \cite{BDI1}. 
Since rotational surfaces of constant
mean curvature are not totally umbilical,
there are no extrinsic spheres
(totally umbilical surfaces with 
constant mean curvature)
in $(G,g[1])$.

On the other hand,
timelike isometric immersion
of $\mathbf{E}^2_1$
into $H^3_1$ with 
parallel second
fundamental form 
are classified in
p.~93, Corollary in \cite{DN}.
See also \cite{Ma}.
\end{Remark}

\begin{Remark} 
Let $c(t)=(x(t),y(t))$ be a curve
in $H^2(1/2)$
parametrised by the arclength
parameter $t$ and 
$M$ the Hopf cylinder over $c$.
Then it is easy to see that 
$\xi$ is tangent to $M$.
Moreover the horizontal lift 
$c^{\prime}(t)^*$ of the tangent vector
field $c^{\prime}(t)$ of $c$ to $G$
 also tangents to $M$. The tangent space
of $M$ at $(x(t),y(t),\theta)$ is spanned by
$c^{\prime}(t)^*$ and $\xi$. 
Denote by $\mathscr{D}^{\perp}$ 
the distribution spanned by $c^{*}(t)$
and put $\mathscr{D}=\{0\}$.
Then we have
$$
TM=\mathscr{D}
\oplus \mathscr{D}^{\perp}\oplus 
\langle \xi \rangle,
\ \
F(\mathscr{D})\subset \mathscr{D},
\ \ F(\mathscr{D}^{\perp})=T^{\perp}M.
$$
Here $\langle \xi \rangle$ is the distribution
spanned by $\xi$.
Thus the Hopf cylinder $M$ is an {\it anti invariant}
submanifold of $G$
in the sense of \cite{YK}.
\end{Remark}

\vspace{0.3cm}

\noindent
{\bf 3.2} \hspace{0.15cm} Next 
we shall recall 
the notion of {\it conoid} introduced by Kokubu.

\begin{Definition} {\rm (\cite{Ko})}
An immersed surface in $(G,g[1])$ of the form:

\begin{equation}
\varphi(u,v)=
\left (
\begin{array}{cc}
1 & x(u) \\ 0 & 1 
\end{array}
\right )
\left (
\begin{array}{cc}
\sqrt{v} & 0 \\ 0 & 1/\sqrt{v} 
\end{array}
\right)
\left (
\begin{array}{cc}
\cos u  & \sin u \\
-\sin u & \cos u 
\end{array}
\right)
\end{equation}
is called a {\it conoid} in $G$. 
\end{Definition}

If we use the metric $g[-1]$, then 
$(x,\theta)=(x(u),u)$ 
is a curve in the double covering
manifold ${\widetilde S}^2_1$
of
$S^2_1$.  
Hence conoids in $(G,g[-1])$ may be
regarded as Hopf cylinders
over curves in $S^2_1$.

Constant mean curvature Hopf
cylinders in $H^3_1$ over curves in $S^2_1$
are classified by
Barros, Ferr{\'a}ndez, Lucas
and Mer{\~o}no.

\begin{Proposition}
{\rm (Classification 
of CMC Hopf cylinders
\cite{BFLM3})}

Let $c$ be a unit speed curve in $S^2_1(1/2)$ 
with curvature $\kappa$ and
$M_c$ the Hopf cylinder over $c$ in $H^3_1$.
Then $M_c$ is
of constant mean curvature if and only if
$c$ is a semi-Riemannian circle in $S^2_1(1/2)$. 
The mean curvature of $M_c$ is
$H=\kappa/2$.
\newline {\rm(1)}
$M_c$ is a minimal complex circle if $\kappa=0$,
\newline {\rm(2)}
$M_c$ is a non-minimal complex circle if $0<\kappa^2<4$,
\newline {\rm(3)}
$M_c$ is a Hopf cylinder over a pseudo-horocycle if $\kappa^2=4$,
\newline {\rm(4)}
$M_c$ is the semi-Riemannian product 
$H^1_1(-r^2)\times S^1_1(r^2-1)$
if $\kappa^2>4$,
\newline {\rm(5)}
$M_c$ is the Riemannian product 
$H^1(-r^2)\times H^1(r^2-1)$ with $r$ such that
$$
\frac{1-2r^2}{r\sqrt{1-r^2}}=\kappa.
$$

\end{Proposition}

\vspace{0.5cm}

On the other hand in
$(G,g[1])$, 
Kokubu obtained the following

\begin{Proposition}{\rm (\cite{Ko})}
The only {\rm(}complete{\rm)} minimal 
conoids in $(G,g[1])$ are
helicoidal surfaces{\rm:}
$$
\varphi(u,v) =
\left (
\begin{array}{cc}
1 & \mu u+a \\ 0 & 1 
\end{array}
\right)
\left(
\begin{array}{cc}
\sqrt{v} & 0 \\ 0 & 1/\sqrt{v} 
\end{array}
\right)
\left(
\begin{array}{cc}
\cos u  & \sin u \\
-\sin u & \cos u 
\end{array}
\right)
$$
$$
=\sigma^\mu_{u}
\left( \
\left(
\begin{array}{cc}
1 & a \\ 0 &1 
\end{array}
\right)
\left(
\begin{array}{cc}
\sqrt{v} & 0 \\ 0 & 1/\sqrt{v} 
\end{array}
\right)
\
\right).
$$

Namely these minimal conoids are $\{\sigma^\mu_t\}$-orbits
of a line $\{(a,y,0) \in H^2 \times 
S^1 \> | \>
y>0\}$.  
In particular $\varphi$ is an imbedding.

\end{Proposition}

The results in this section 
motivate us to study the class of
surfaces which will be
introduced 
in the next section.

\section{Surfaces derived from curves in the lightcone}

In this section,
we shall introduce a new class of surfaces in $G$.
As we saw before,
$\mathrm{Ad}$-orbits of vectors
in ${\mathfrak g}
={\mathfrak s}{\mathfrak l}_{2}{\mathbf R}$
are classified in three types.
The $\mathrm{Ad}$-orbit
of a spacelike [resp. timelike]
vector is a hyperbolic $2$-space
[resp. Lorentz sphere].
The $\mathrm{Ad}$-orbit
of a null vector
is the lightcone.
In the preceding section, 
we saw that two kinds of surfaces,
``rotational surfaces" and 
``conoids" coincide Hopf cylinders over
curves in hyperbolic $2$-space or
Lorentz sphere.
It seems to be interesting to study
surfaces obtained by curves in
$\mathrm{Ad}$-orbit of a null vector,
{\it i.e.,} the lightcone.
This section is devoted to study
such surfaces.
\vspace{0.2cm}

Let $c$ be a curve in lightcone $\Lambda$.
Then its inverse image 
$M$ in $H^3_1=(G,g[-1])$ is given by
\begin{equation}\label{null}
\varphi(u,v)=
\left(
\begin{array}{cc}
1 & v \\ 0 &1 
\end{array}
\right)
\left(
\begin{array}{cc}
\sqrt{y(u)} & 0 \\
0 & 1/\sqrt{y(u)} 
\end{array}
\right)
\left(
\begin{array}{cc}
\cos u & \sin u \\
-\sin u & \cos u
\end{array}
\right).
\end{equation}
The partial derivatives of $\varphi$ are
$$
\varphi_{*}
\frac{\partial}{\partial u}=
\frac{y^{\prime}}{2y}\epsilon_2+
\epsilon_3,\ \ 
\varphi_{*}
\frac{\partial}{\partial v}=
\frac{1}{2y}
\left (\epsilon_1+
\epsilon_3
\right ).
$$
The induced metric 
$\mathrm{I}[\nu]$ of $M$ is
$$
\mathrm{I}[\nu]=
\left
\{
\nu+\left(
\frac{y^{\prime}(u)}
{2y(u)}
\right)^2
\right \}
du^2+
\frac{\nu}{y(u)}dudv+
\frac{1+\nu}{4y(u)^2}dv^2.
$$
The determinant of 
$\mathrm{I}[\nu]$ is
$$
\det 
\mathrm{I}[\nu]=
\frac{1}{16y(u)^4}
\left \{
(1+\nu)y^{\prime}(u)^2+4\nu
y(u)^2
\right\}.
$$
In particular
$\det \mathrm{I}[-1]=-1/(4y^2)$,
hence $(M,\varphi)$ is timelike in $H^3_1$.
Direct computations using (\ref{con}) show that
$$
\nabla_{\partial_u}\varphi_{*}
\frac{\partial}{\partial u}=
-\nu
\left(
\frac{y^\prime}
{y}
\right)
\epsilon_1
+
\left(
\frac{y^\prime}
{2y}
\right)^{\prime}
\epsilon_2,
$$
$$
\nabla_{\partial_u}\varphi_{*}
\frac{\partial}{\partial v}=
\frac{1}{4y^2}
\left \{
-y^\prime(\nu+2)\epsilon_1
+2\nu y\epsilon_2-y^\prime
\epsilon_3
\right \},
$$
$$
\nabla_{\partial_v}\varphi_{*}
\frac{\partial}{\partial v}=
\frac{\nu+1}{2y^2}
\epsilon_2.
$$

The unit normal vector field 
${\mathbf n}[\nu]$ is
$$
\mathbf{n}[\nu]=\frac{1}
{\sqrt{1+(1+\nu)(\frac{y^\prime}{2y})^2}}
\left(
\frac{y^\prime}{2y}\epsilon_1
+\epsilon_2-
\frac{\nu y^\prime}{2y}
\epsilon_3
\right ).
$$

Let us denote by
${\mathrm I}\!{\mathrm I}={\mathrm I}\!{\mathrm I}[\nu]$ the
{\it second fundamental 
form} derived form 
$\mathbf{n}[\nu]$.
The second fundamental form
${\mathrm I}\!{\mathrm I}$ is defined by the
{\it Gau{\ss} formula}: 
$$
\nabla_{X}\varphi_{*}Y=\varphi_{*}
(\nabla^M_{X}Y)+
{\mathrm I}\!{\mathrm I}(X,Y)\mathbf{n},\ \
X,Y \in \mathfrak{X}(M).
\eqno{\rm (G)}
$$
Here $\nabla^M$ is the 
Levi-Civita connection of
$(M,\mathrm{I}[\nu])$.

Put $\alpha=
\sqrt{1+(1+\nu)\{y^{\prime}/(2y)\}^2}$.
Then $\det \mathrm{I}[\nu]=
\nu\alpha^{2}/(4y^2)$.

The second fundamental
form ${\mathrm I}\!{\mathrm I}$ is described by the
following formulae:
$$
{\mathrm I}\!{\mathrm I}(\frac{\partial}{\partial u},
\frac{\partial}{\partial u})=
\frac{-(1+\nu)
y^{\prime}(u)^2+y^{\prime \prime}(u)y(u)}
{2 \alpha y(u)^2},
$$
$$
{\mathrm I}\!{\mathrm I}(\frac{\partial}{\partial u},
\frac{\partial}{\partial v})
=\frac
{-(1+\nu)y^{\prime}(u)^2+4\nu y(u)^2}
{8\alpha y(u)^3},
$$
$$
{\mathrm I}\!{\mathrm I}
(\frac{\partial}{\partial v},
\frac{\partial}{\partial v})
=\frac{(1+\nu)}{2 \alpha y(u)^2}.
$$
The mean curvature $H[\nu]$ 
of $\varphi$ is 
\begin{equation}\label{Mean}
H[\nu]=
\frac{1}
{4 \alpha^3 y(u)^2}
\left \{
(1+\nu)
y^{\prime \prime}(u)y(u)+
4y(u)^2
\right \}.
\end{equation}
Here we used the formula:
$$
H[\nu]=\frac{1}{2}\>\mathrm{tr}
\{{\mathrm I}\!{\mathrm I}[\nu] \cdot 
\mathrm{I}[\nu]^{-1}
\}.
$$

\noindent
{\bf Case 1}: $\nu=1$ 

From (\ref{Mean}),
we have $\varphi$ is
minimal if and only if
$$
y^{\prime \prime}=-2y.
$$

\begin{Theorem}
Let $\varphi(u,v)$ an 
immersed surface in $(G,g[1])$ 
obtained by taking inverse image 
of a curve in $\Lambda$
which is parametrised as 
{\rm (}{\rm \ref{null}}{\rm )}.
Then $\varphi$ is minimal if and only if
$\varphi$ is the inverse image of
$$
(A\cos (\sqrt{2}u)+B\sin (\sqrt{2}u),u)
\in \mathbf{R}^{+}\times S^1.
$$
\end{Theorem}

\noindent
{\bf  Case 2:} $\nu=-1$

On the other hand, in $(G,g[-1])$, 
$\varphi$ has constant mean curvature $1$
and Gau{\ss}ian curvature $0$.
Denote by $\mathcal{D}$ the discriminant
of the characteristic equation:
$$
\det (tI-S)=0
$$
for the shape operartor 
$S={\mathrm I}\!{\mathrm I} 
\cdot \mathrm{I}^{-1}$.
Then $\mathcal{D}$ is given by the 
following formula:
$$
\mathcal{D}=H^2-K-1.
$$
Thus $(M,\varphi)$ has real and repeated 
principal curvatures in $H^3_1$.
Hence 
$\varphi$ is 
a $B$-scroll of a null Frenet curve with
constant torsion $1$ in $H^3_1$. 
In particular $M$ is flat
totally umbilical timelike surface
if and only if it is a $B$-scroll
of a null geodesic with constant torsion $1$.  
(See Theorem 3 in \cite{DN}. )

Comparing the first and second fundamental
forms we have the following

\begin{Proposition}
Let $\varphi(u,v)$ an 
immersed surface in $H^3_1$ 
obtained by taking inverse image of a curve in $\Lambda$
which is parametrised as {\rm (}{\rm \ref{null}}{\rm )}.
Then $\varphi$ is totally umbilical if and only if
$y$ is a solution to
\begin{equation}\label{ODE}
y^{\prime \prime}-
\frac{(y^\prime)^2}{2y}
+2y=0.
\end{equation}
\end{Proposition}

The ordinary differential equation 
(\ref{ODE}) with $y>0$ can be solved explicitly.
In fact let us introduce an auxiliary function
$\mathscr{T}$ by 
$$
\mathscr{T}(u):=\frac{d}{du}\log y(u).
$$
Then (\ref{ODE}) is rewritten as
$$
\mathscr{T}^{\prime}
+\frac{1}{2}\mathscr{T}^{2}+2=0.
$$
The general solutions of 
this ordinary equation are given 
explicitly by
$$
\mathscr{T}(u)=-2\tan (u+u_0),\ u_0\in \mathbf{R}.
$$
Thus the solutions $y$ to (\ref{ODE}) are
given by
$$
y(u)=A\> 
\cos^{2} 
(u+u_0),\ A>0.
$$

\begin{Theorem}\label{umbilicScroll}
Let $\varphi(u,v)$ an immersed surface in 
$H^3_1$ obtained by taking inverse image of 
a curve in $\Lambda$.
Then $\varphi$ is a 
$B$-scroll of a null Frenet curve 
with constant torsion $1$ in $H^3_1$. 
In particular $\varphi$ is totally umbilical
if and only if $\varphi$ is the inverse image of
the curve{\rm :}
$$
(A\> 
\cos^{2} 
(u+u_0),u)
\in \mathbf{R}^{+} \times S^1, \ \ A>0.
$$ 
\end{Theorem}

\begin{Remark}(Weierstra{\ss}-type 
representations for surfaces in $H^3_1$)
\begin{enumerate}
\item
Hong \cite{Ho} obtained a Bryant-type
representation formula
\newline
\noindent
for
timelike constant mean curvature $1$ surfaces
in $H^3_1$.
\item
Balan and 
Dorfmeister \cite{BD} established a loop group theoretic
\newline
\noindent
Weierstra{\ss}-type 
representation (so-called {\it DPW representation})
for harmonic maps of Riemann surface into 
general Lie group with biinvariant 
\newline
\noindent
semi-Riemannian metric. 
Their general scheme is 
applicable to maximal (spacelike) surfaces in 
$H^3_1=(\mathrm{SL}_{2}{\mathbf R},g[-1])$.  
\end{enumerate}
\end{Remark} 

\begin{Remark}
Hopf cylinders over curves in $H^2$
[resp. $S^2_1$] are surfaces in $G$ which are
invariant under $K$-action [resp. $A\mathbf{Z}_2$-action].
Surfaces considered in this section are
invariant under $N$-action.
Thus all the surfaces investigated in preceding section
and present section are invariant under
$1$-dimensional closed subgroup of the isometry
group $G\times K$.
In \cite{FMP}, Figueroa, Mercuri and
Pedrosa classified all
constant mean curvature surfaces
in the Heisenberg group which are invariant
under $1$-dimensional closed subgroups
of the isometry group.
Some results in 
\cite{FMP} are independently
obtained in
\cite{IKOS2}.
Recently S.~D.~Pauls studied minimal surfaces 
in the Heisenberg group with
Carnot-Carath\'eodory metric \cite{Pau} 

\end{Remark}

\section{Tangential Gau{\ss} maps}
\noindent
{\bf 5.1} \hspace{0.15cm} Let 
$(N^n,g_N)$ be a Riemannian $n$-manifold
and $O(N)$ the orthonormal frame bundle of
$N$. As is well known,
$O(N)$ is a principal $\mathrm{O}(n)$-bundle
over $N$.

Denote by $Gr_{\ell}(T_pN)$ be the Grassmannian manifold
of $\ell$-planes in the tangent space $T_pN$ of $N$ at $p
\in N$.
The set $Gr_{\ell}(TN):=\cup_{p \in N}
Gr_{\ell}(T_pN)$ of all $\ell$-planes 
in the tangent bundle
$TN$ admits a structure of fibre bundle
over $N$. 
In fact, 
$Gr_{\ell}(TN)$ is a fibre bundle associated
to $O(N)$:
$$
Gr_{\ell}(TN)=O(N) \times_{{\mathrm O}(n)}Gr_{\ell}(\mathbf{E}^n)
$$
whose standard fibre 
is the Grassmannian manifold
$Gr_{\ell}(\mathbf{E}^n)$
of $\ell$-planes in Euclidean
$n$-space.
This fibre bundle 
$Gr_{\ell}(TN)$ is called the
{\it Grassmannian bundle}
of $\ell$-planes over $N$.

The canonical 1-form of $O(N)$ and 
the Levi-Civita connection 1-forms of $g_N$
naturally induces an 
invariant Riemannian metric
$\langle \cdot,\cdot
\rangle$ on $Gr_{\ell}(TN)$.
with respect to this metric the 
projection 
$pr:Gr_{\ell}(TN) \to N$ becomes a
Riemannian submersion with 
totally geodesic fibres.
For more details about the 
metric, see Jensen and Rigoli
\cite{JR} and Sanini \cite{Sa}.

\begin{Definition}\label{Gaussmapdef}
Let $\varphi:M^m\to N^n$ be an immersed submanifold.
Then the (tangential) {\it Gau{\ss} map} 
$\psi:M \to Gr_{m}(TN)$ is defined by
$$
\psi(p):=\varphi_{*p}(T_{p}M) \in Gr_{m}(T_{p}N),\ \ p \in M.
$$ 
\end{Definition}

\begin{Remark}
In case the ambient Riemannian $n$-manifold $N$ is
a Lie group with left invariant metric
and $M$ is a hypersurface,
we can introduce another kind
of Gau{\ss} map.

Let $G$ be an $n$-dimensional Lie group
with left invariant metric.
For an immersed hypersurface $\varphi:M \to G$
with unit normal ${\mathbf n}$,
the {\it normal}
{\it Gau{\ss}} {\it map} 
$\Upsilon$ of $M$ is
a smooth map into the unit $(n-1)$-sphere in 
the Lie algebra ${\mathfrak g}$ of $G$ defined by
$$
\Upsilon(p):=L_{\varphi(p)*}^{-1}{\mathbf n}_{p} \in 
S^{n-1} \subset \mathfrak{g}.
$$ 
In our study for surfaces in $(G,g[1])$,
to distinguish the Gau{\ss} maps into the
$Gr_{2}(TG)$ from the normal Gau{\ss} maps,
we use the name ``tangential Gau{\ss} maps" for
the Gau{\ss} maps defined in Definition \ref{Gaussmapdef}. 
\end{Remark}

\vspace{0.3cm}

\noindent
{\bf 5.2}
\hspace{0.15cm} 
Here we recall and colllect
fundamental ingredients in 
the theory of {\it harmonic maps} 
from the lecture note \cite{EL}
by Eells and Lemaire.

Let $(M,g_M)$ and $(P,g_P)$ be Riemannian manifolds.
And let $f:M \to P$ be a smooth map of a manifold $M$
into $P$. The {\it energy density}
$e(f)$ of $f$ is a smooth function
on $M$ defined by
$e(f):=|df|^2/2$. It is obvious that
$e(f)=0$ if and only if $f$ is constant.

The {\it energy} $E(f)$ of $f$
is 
$$
E(f):=\int_{M} e(f) \> dV_{M}.
$$
Here $dV_M$ is the volume element of $(M,g_M)$. 

The {\it tension field} 
$\tau(f)$ of $f$ is a smooth section of 
$f^{*}(TP)$ defined by
$$
\tau(f):=\mathrm{tr}\> \nabla d f.
$$
It is known that $f$ is a critical point of 
the energy if and only if $\tau(f)=0$.

A map $f$ is said to be a
{\it harmonic map} if $\tau(f)=0$.

Baird and Eells introduced
the notion of {\it stress-energy tensor} 
in \cite{BE}.
The stress-energy tensor 
$\mathcal{S}(f)$
of a map $f$ is a symmetric 
(0,2)-tensor field on $M$ defined by 
$$
\mathcal{S}(f):=e(f)g_{M}-f^{*}g_{P}.
$$

In particular in case $\dim M=2$ and $f$ is nonconstant,
$f$ is conformal if and only if
$\mathcal{S}(f)=0$.

Since $\mathcal{S}(f)$ is symmetric
$(0,2)$-tensor field, 
the divergence  
$\mathrm{div} \> \mathcal{S}(f)$ of 
$\mathcal{S}(f)$ can be defined
by the formula: 
$$
\mathrm{div} \> \mathcal{S}(f):=
\mathbb{C}_{13}(\nabla \mathcal{S}(f)\>).
$$
(See p.~86 in \cite{O})

Here $\mathbb{C}_{13}$ is the 
metric contraction operator in the
1st and 3rd entries. See p.~83 in 
\cite{O}.
The divergence of $\mathcal{S}(f)$ is 
given explicitly by \cite{BE}:
$$
\mathrm{div} \> \mathcal{S}(f)=-g_{P}(\tau(f),d f).
$$
Thus if $f$ is a harmonic map then its stress-energy tensor is
conservative.

\vspace{0.2cm}

\noindent
{\bf 5.3} 
\hspace{0.15cm} Next 
we recall the notion 
of {\it vertically harmonic map}
\cite{Wc}.

Let $(P,g_P)$ be a Riemannian manifold
and $pr:(P,g_P) \to (N,g_N)$ a Riemannian submersion.
With respect to the metric $g_P$,
the tangent bundle $TP$ of $P$ is decomposed as:
$$
T_{u}P=\mathcal{H}_{u}\oplus \mathcal{V}_u,\ u \in P.
$$
Here $\mathcal{V}_u:=\mathrm{Ker}\ (pr_{*})_{u}$
and $\mathcal{H}_u=\mathcal{V}_{u}^{\perp}$ are called
the {\it vertical subspace} and
{\it horizontal subspace} of $T_{u}P$ at $u$ 
respectively.

Now let $f:(M,g_M) \to (P,g_P)$ be a smooth map.
With respect to the Riemannian 
submersion $pr$, 
$\tau(f)$ is decomposed into its
horizontal and vertical
components:
$$
\tau(f)=\tau^{\mathcal H}(f)+
\tau^{\mathcal V}(f).
$$
The map $f$ is said to be a
{\it vertically harmonic map}
if the vertical component
$\tau^{\mathcal V}(f)$ vanishes.

In case $f:M=N \to P$ is a section of $P$,
{\it i.e}., a smooth map satisfying $pr \circ f=\mathrm{identity}$,
C.~M.~Wood \cite{Wc} showed that 
the vertical harmonicity for maps
is equivalent to the criticality
for the vertical energy under the vertical
variations.

\vspace{0.2cm}

\noindent
{\bf 5.4} \hspace{0.15cm} Now
we investigate
harmonicity of
tangential Gau{\ss} maps for surfaces in 
$(G,g[1])$. 

The following fundamental
result is  
due to Sanini
(See (3.2)--(3.3) in \cite{Sa}).

\begin{Lemma}\label{SaniniLemma}
Let $N$ be a Riemannian $3$-manifold
and $\varphi:M \to N$ an immersed surface
with unit normal vector field
${\mathbf n}$.
Take a principal frame field
\newline
\noindent
$\{e_1,e_2,e_3=\mathbf{n}\}$,
{\rm i.e.}, an orthonormal
frame field such that 
$\{e_1,e_2 \}$ diagonalise
the shape operator. Put
$$
R_{ijkl}=
g_{N}(
R(e_i,e_j)e_k,e_l)
$$
and denote by $\psi$ the tangential
Gau{\ss} map of $(M,\varphi)$.
Then the following holds.

\vspace{0.2cm}

{\rm (1) }
The tangential Gau{\ss} map
$\psi$ 
is conformal if and only if 
$(M,\varphi)$ is totally umbilical or minimal.

\vspace{0.2cm}

{\rm (2)} 
Assume that $(M,\varphi)$ has
constant mean curvature. Then 
$\psi$ is vertically harmonic 
if and only if
$R_{1213}=R_{2123}=0$.
Moreover when $(M,\varphi)$ is minimal,
$\psi$ 
is harmonic if and only if,
in addition, 
$R_{3113}=R_{3223}=0$.

\vspace{0.2cm}

{\rm (3)} 
Assume that the mean curvature is nonzero
constant. Then 
the tangential Gau{\ss} map
is vertically harmonic if
and only if the stress energy
tensor $\mathcal{S}(\psi)$ of the tangential
Gau{\ss} map $\psi$ is conservative 
{\rm (}divergence free{\rm )}.
\end{Lemma}

Sanini applied this Lemma to 
surfaces in 3-dimensional Heisenberg group
with canonical left invariant
metric \cite{Sa}.

Lemma \ref{SaniniLemma} together with the nonexistence
of extrinsic spheres 
(See Remark \ref{parallel} and 
\cite{BDI1}) implies the
following.

\begin{Corollary}
Let $M$ be a constant mean curvature
surface in $(G,g[1])$. Then $M$ is minimal if and
only if its tangential Gau{\ss} map is conformal.
\end{Corollary} 

The following is the 
main result of this section\footnote{This result is generalised to 
$3$-dimensional Sasakian space forms by M.~Tamura (Comment. Math. Univ. St. Pauli \textbf{52} (2003), no.~2, 117--123.}.

\begin{Theorem}
Let $M$ be a surface in $(G,g[1])$
with constant mean curvature.
Then the tangential Gau{\ss} map of $M$ is
vertically harmonic if and only if
$M$ is a Hopf cylinder {\rm (}
rotational surface{\rm )} of constant
mean curvature.
Hopf cylinders with nonzero 
constant mean curavture are 
{\rm(}only{\rm)} 
constant mean curvature surfaces
whose tangential 
Gau{\ss} map are vertically harmonic but
nonharmonic and have 
conservative stress-energies.

In particular the only minimal surface
in $(G,g[1])$ with vertically harmonic
tangential Gau{\ss} map is a
Hopf cylinder over a geodesic.
In this case the tangential
Gau{\ss} map is a harmonic map.
\end{Theorem}
{\bf Proof.}
Let $\varphi:M \to (G,g[1])$ 
be a surface with
constant mean curvature
and unit normal vector
field ${\mathbf n}$.
Denote by $\theta^3$ the dual
one-form of ${\mathbf n}$.
Express
$\theta^3$ by
$$
\theta^3=a\> \omega^1+
b\> \omega^2+c\> \omega^3,
\ \
a^2+b^2+c^2=1
$$
in terms of the coframe field 
(\ref{cobase}).
\begin{enumerate}

\item{Case 1 $c\ne 0$:}
In this case,
$$
v_1=-c\epsilon_2+b\epsilon_3, \
v_2=(b^2+c^2)\epsilon_1-ab\epsilon_2
-ac\epsilon_3
$$
gives a orthogonal frame
field of $M$.

Direct computations show the following
formulae:
$$
g[1](R(v_1,v_2)v_1,\mathbf{n})=
8ac^{2}(b^2+c^2),\ \
g[1](
R(v_1,v_2)v_2,
\mathbf{n})=
8bc(b^2+c^2).
$$
Take a principal frame 
$\{e_1,e_2\}$. Then 
$\{e_1,e_2\}$ is expressed as
\begin{equation}\label{VE}
e_1=\cos \mu \frac{v_1}{|v_1|}
+\sin \mu \frac{v_2}{|v_2|},\ \
e_2=-\sin \mu
\frac{v_1}{|v_1|}
+\cos \mu \frac{v_2}{|v_2|}.
\end{equation}
Then we have
$$
R_{1213}=
\frac{8c(b^2+c^2)}
{|v_1||v_2|}
\left(
\frac{ac}{|v_1|}\cos \mu+
\frac{b}{|v_2|}\sin \mu
\right),
$$

$$
R_{2123}=
\frac{8c(b^2+c^2)}
{|v_1||v_2|}
\left(
\frac{ac}{|v_1|}\sin \mu-
\frac{b}{|v_2|}\cos \mu
\right).
$$
From these we have $\tau^{\mathcal V}(\psi)=0$
if and only if $a=b=0$.
Hence $\theta^3=-\eta$.
Namley $M$ is an integral surface
of the distribution $\eta=0$,
but this is impossible, 
since
$\eta$ is contact.
(See p.~36, Theorem
in \cite{Bl}).

\vspace{0.2cm}

\item{Case 2 $c=0$:}
Since $a^2+b^2=1$, we may write
$a=\cos \phi,\ b=\sin \phi$.

In this case $u_1=
\sin \phi\epsilon_1-\cos \phi \epsilon _2,\ 
u_2=\epsilon_3$
are orthonormal and tangent to $M$.
The unit normal ${\mathbf n}$ is given by
${\mathbf n}=\cos \phi \epsilon_1+
\sin \phi \epsilon_2$.
Then we have
\begin{equation}\label{uR}
R(u_1,u_2)u_1=-\sin^{2} \phi \epsilon_3,\ \ 
R(u_2,u_1)u_2=-\sin \phi \epsilon_1
+\cos \phi \epsilon_2.
\end{equation}

Let us denote by $\mu$ the angle
between the principal frame 
$\{e_1,e_2\}$ and $\{u_1, u_2 \}$,
{\it i.e.,} 
\begin{equation}\label{UE}
e_1=\cos \mu \>u_1+\sin \mu \>u_2,\ \
e_2=-\sin \mu \>u_1+\cos \mu \>u_2.
\end{equation}
Using (\ref{uR}) and (\ref{UE}), 
we have 
$R_{1213}=R_{2123}=0$.
Thus $\tau^{\mathcal V}(\psi)=0$
is fulfilled automatically for
$M$ with $c=0$.

We have shown in 
\cite{BDI1} that
constant mean curvature surfaces with 
$c=0$ are Hopf cylinder of constant
mean curvature.
See the proof of Thereom in \cite{BDI1}.

Furthermore,
the second fundamental form
${\mathrm I}\!{\mathrm I}$ of $M$ relative to ${\mathbf n}$ is
given by ({\it cf}. (5) and (8) in \cite{BDI1}) 
\begin{equation}\label{sff}
{\mathrm I}\!{\mathrm I}(u_1,u_1)=2H,\ 
{\mathrm I}\!{\mathrm I}(u_1,u_2)=1,\ \
{\mathrm I}\!{\mathrm I}(u_2,u_2)=0.
\end{equation}

Next we see the case $\psi$ 
is harmonic. 
Using (\ref{uR}) and (\ref{UE}) 
again, we have 
$$
R_{3113}=-7 \cos^{2}\mu+\sin^{2}\mu,\ \
R_{3223}=-7 \sin^{2}\mu+\cos^{2}\mu.
$$
Thus $R_{3113}=R_{3223}$ if only if
$\mu=\pm \pi/4$. Without
loss of generality, we may assume
$\mu=\pi/4$.
In this case the principal
frame $\{e_1,e_2\}$ is given by
$$
e_1=\frac{1}{\sqrt{2}}(u_1+u_2),\ \
e_2=\frac{1}{\sqrt{2}}(-u_1+u_2).
$$
By definition, ${\mathrm I}\!{\mathrm I}(e_1,e_2)=0$.
On the other hand, direct computation using
(\ref{sff}) shows 
${\mathrm I}\!{\mathrm I}(e_1,e_2)=-H$. Thus 
a constant mean curvature
surface $M$ with $c=0$
satisfying
$R_{3113}=R_{3223}$ is minimal.
\end{enumerate}
Conversely one can check that every rotational surface of
constant mean curvature has vertically
harmonic tangential Gau{\ss} map
and when $H\not=0$, 
the tension field does not vanish
by direct computations. 
It is also straightforward to
check that every minimal Hopf
cylinder has harmonic
tangentail Gau{\ss} map.
$\Box$

\vspace{1.0cm}

\noindent
J.~Inoguchi 
\newline
\noindent
Department of Applied Mathematics
\newline
\noindent
Fukuoka University
\newline
\noindent
Nanakuma, 
Fukuoka, 814-0180 
\newline
\noindent
Japan

\vspace{0.15cm}
\noindent
{\tt inoguchi@bach.sm.fukuoka-u.ac.jp}

\vspace{1.0cm}

\noindent
Current address
\newline
\noindent
Department of Mathematical Sciences
\newline
\noindent
Yamagata University
\newline
\noindent
Yamagata 990-8560
\newline
\noindent
Japan

\vspace{0.15cm}
\noindent
{\tt inoguchi@sci.kj.yamagata-u.ac.jp}


\begin{thebibliography}{99}

\bibitem{AFL}
L.~J.~Al{\' \i}as,
A.~Ferr{\'a}ndez
and P.~Lucas,
$2$-type surfaces in
$S^3_1$ and $H^3_1$,
Tokyo J. Math.
{\bf 17} (1994),
447--454.

\bibitem{BE}
P.~Baird and J.~Eells,
A conservation law
for harmonic maps,
in:
{\sl
Geometry Symposium, Utrecht} 1980 (Utrecht, 1980),
Lecture Notes in Math.
{\bf 894} (1981),
Springer Verlag,
pp.~1--25.


\bibitem{BD}
V.~Balan and J.~Dorfmeister,
A Weierstrass-type
representation for harmonic 
maps from Riemann surfaces
to general Lie groups,
Balkan J. Geom. Appl.
{\bf 5} (2000), 7--37,
(http://www.emis.de/journals/BJGA/5.1/2.html). 

\bibitem{BFLM}
M.~Barros, A.~Ferr{\'a}ndez,
P.~Lucas and M.~Mero{\~n}o, 
{\rm Hopf cylinders,
$B$-scrolls and solitons of
the Betchov-da Rios equation
in the $3$-dimensional anti-de Sitter
space},
C. R. Acad. Sci.
Paris S{\'e}rie I
{\bf 321} (1995),
505--509.


\bibitem{BFLM3}
M.~Barros, A.~Ferr{\'a}ndez,
P.~Lucas and M.~Mero{\~n}o, 
{\rm Solutions of the Betchov-da Rios
soliton equation in the anti de Sitter 
$3$-space}, in:
{\sl New Approaches in 
Nonlinear Analysis},
(Th.~M.~Rassias ed.), Hadronic
Press, Inc., Palm Harbor, Florida, 1999, 
pp.~51--71.


\bibitem{BDI1}
M.~Belkhelfa, F.~Dillen and
J.~Inoguchi,
Parallel surfaces in the 
real special linear
group $SL(2,\mathbb{R})$,
Bull. Austral. Math.
Soc. {\bf 65} (2002),
183--189.

\bibitem{Bl}
D.~E.~Blair,
{\sl Contact Manifolds in 
Riemannian Geometry},
Lecture Notes in Math.
{\bf 509} (1976),
Springer Verlag,
Berlin.

\bibitem{DN}
M.~Dajczer and K.~Nomizu,
{\rm On flat surfaces in $S^3_1$ and $H^3_1$},
in:
{\sl Manifolds and Lie Groups--papers
in honor of Yozo Matsushima}
(J.~Hano et al eds.),
Progress in Math.
{\bf 14} (1981),
Birkh{\"a}user, Boston,
pp.~71-108.


\bibitem{EL}
J.~Eells and L.~Lemaire,
{\sl Selected Topics in Harmonic Maps},
Regional Conference Series in Math.
{\bf 50} (1983),
Amer. Math. Soc.
 
\bibitem{F}
A.~Ferr{\'a}ndez,
Riemannian versus Lorentzian submanifolds,
some open problems,
in:
{\sl Proc. Workshop on
Recent Topics in Differential Geometry,
Santiago de Compostela},
Depto. Geom.
y Topolog{\' \i}a,
Univ. Santiago de
Compostela,
{\bf 89} (1998),
pp.~109--130.


\bibitem{FMP}
C.~B.~Figueroa,
F.~Mercuri
and
R.~H.~L.~Pedrosa,
Invariant minimal surfaces
of the Heisenberg groups,
Ann. Mat. Pura Appl.
{\bf 177} (1999),
173--194.

\bibitem{Grod}
C.~Gorodski,
Delaunay-type surfaces in the
$2 \times 2$ real unimodular group,
Ann. Mat. Pura Appl.
{\bf 180} (2001),
211--221.

\bibitem{Ho}
J.~Q.~Hong, Timelike surfaces with mean curvature one in anti 
de Sitter $3$-space,
K{\=o}dai Math. J.
{\bf 17} (1994), 341-350.



\bibitem{IKOS2}
J.~Inoguchi, T.~Kumamoto,
N.~Ohsugi and
Y.~Suyama,
Differential geometry of curves and surfaces in 
$3$-dimensional homogeneous spaces 
$\mathrm{I}\!\mathrm{I}$,
Fukuoka Univ. Sci. Rep.
{\bf 30}
(2000),
17--47. 


\bibitem{IKOS3}
J.~Inoguchi, T.~Kumamoto,
N.~Ohsugi and
Y.~Suyama,
{\rm Differential geometry of curves and surfaces in 
$3$-dimensional homogeneous spaces 
$\mathrm{I}\!\mathrm{I}\!\mathrm{I}$},
Fukuoka Univ. Sci. Rep.
{\bf 30} (2000),
131--160.

\bibitem{IKOS4}
J.~Inoguchi, T.~Kumamoto,
N.~Ohsugi and
Y.~Suyama,
{\rm Differential geometry of curves and surfaces in 
$3$-dimensional homogeneous spaces 
$\mathrm{I}\!\mathrm{V}$},
Fukuoka Univ. Sci. Rep.
{\bf 30} (2000),
161--168.

\bibitem{JR}
G.~Jensen and M.~Rigoli,
Harmonic Gauss maps,
Pacific J. Math.
{\bf 136} (1989),
261--282.



\bibitem{Ko}
M.~Kokubu,
{\rm On minimal surfaces in the real
special linear group}
$SL(2,{\mathbf R})$,
Tokyo J. Math.
{\bf 20} (1997), 
287--297.

\bibitem{Ma}
M.~Magid,
{\rm Isometric immersions of
Lorentz space with parallel second fundamental forms},
Tsukuba J. Math.
{\bf 8}
(1984),
31--54.


\bibitem{O}
B.~O'Neill,
{\sl Semi-Riemannian Geometry with Application
to Relativity},
Academic Press, 1983.

\bibitem{Pau}
S.~D.~Pauls,
Minimal surfaces in the Heisenberg
group,
Geom. Dedicata \textbf{104} (2004), 201--231, 
math.DG/0108048.

\bibitem{P}
U.~Pinkall,
Hopf tori in $S^3$,
Invent. Math.
{\bf 81} (1985),
379--386.


\bibitem{Sa}
A.~Sanini,
Gauss maps of a surface
of the Heisenberg group,
Boll. Un. Mat. Ital.
(7) {\bf 11}-B (1997),
suppl. fasc.
2, 79--93.

\bibitem{Tak}
T.~Takahashi,
Sasakian manifold with 
pseudo-Riemannian metric,
T\^ohoku Math. J. (2)
{\bf 21} (1969),
271--290.


\bibitem{TV}
F.~Tricerri and L.~Vanhecke,
{\sl Homogeneous Strctures
on Riemannian manifolds},
London Math. Soc.
Lecture Note Series
{\bf 83},
 Cambridge University Press,
 Cambridge,
1983.

\bibitem{Wc}
C.~M.~Wood,
The Gauss section of 
a Riemannian immersion,
J. London Math. Soc. (2),
{\bf 33}
(1986),
157--168.


\bibitem{YK}
K.~Yano and M.~Kon,
{\sl $CR$-submanifolds of Kaehlerian and
Sasakian Manifolds},
Birkh\"auser, Boston,
1983


\end{thebibliography}
\end{document}